\documentclass[11pt]{article}
\usepackage{amssymb}
\usepackage{graphicx}
\usepackage{amsmath}
\usepackage{amssymb}
\usepackage{graphicx}
\usepackage{amsmath}

\newcommand{ \Spec}{\mbox{\textup{Spec }}}

\newcommand{ \Ass}{\mbox{\textup{Ass}}}
\newcommand{ \Rad}{\mbox{\textup{Rad}}}

\newcommand{ \Rees}{\mbox{\textup{Rees }}}

\newcommand{ \altitude}{\mbox{\textup{altitude}}}

\newenvironment{proof}[1][Proof]{\noindent{\textbf{#1.}}
}{\ \rule{0.5em}{0.5em}}

\def\inbar{\,\vrule height1.5ex width.4pt depth0pt}
\def\inbar{\,\vrule height1.5ex width.4pt depth0pt}
\def\IB{\relax{\rm I\kern-.18em B}}
\def\IC{\relax\hbox{$\inbar\kern-.3em{\rm C}$}}
\def\ID{\relax{\rm I\kern-.18em D}}
\def\IE{\relax{\rm I\kern-.18em E}}
\def\IF{\relax{\rm I\kern-.18em F}}
\def\IG{\relax\hbox{$\inbar\kern-.3em{\rm G}$}}
\def\IH{\relax{\rm I\kern-.18em H}}
\def\II{\relax{\rm I\kern-.18em I}}
\def\IK{\relax{\rm I\kern-.18em K}}
\def\IL{\relax{\rm I\kern-.18em L}}
\def\IM{\relax{\rm I\kern-.18em M}}
\def\IN{\relax{\rm I\kern-.18em N}}
\def\IO{\relax\hbox{$\inbar\kern-.3em{\rm O}$}}
\def\IP{\relax{\rm I\kern-.18em P}}
\def\IQ{\relax\hbox{$\inbar\kern-.3em{\rm Q}$}}
\def\IR{\relax{\rm I\kern-.18em R}}
\font\cmss=cmss10 \font\cmsss=cmss10 at 7pt
\def\IZ{\relax\ifmmode\mathchoice
{\hbox{\cmss Z\kern-.4em Z}}{\hbox{\cmss Z\kern-.4em Z}}
{\lower.9pt\hbox{\cmsss Z\kern-.4em Z}}
{\lower1.2pt\hbox{\cmsss
Z\kern-.4em Z}}\else{\cmss Z\kern-.4em Z}\fi}
\def\IGa{\relax\hbox{${\rm I}\kern-.18em\Gamma$}}
\def\IPi{\relax\hbox{${\rm I}\kern-.18em\Pi$}}
\def\ITh{\relax\hbox{$\inbar\kern-.3em\Theta$}}
\def\IOm{\relax\hbox{$\inbar\kern-3.00pt\Omega$}}

\begin{document}

\baselineskip 20pt
\pagenumbering{arabic}
\pagestyle{plain}

\newtheorem{defi}{Definition}[section]
\newtheorem{theo}[defi]{Theorem}
\newtheorem{lemm}[defi]{Lemma}
\newtheorem{prop}[defi]{Proposition}
\newtheorem{note}[defi]{Note}
\newtheorem{nota}[defi]{Notation}
\newtheorem{exam}[defi]{Example}
\newtheorem{coro}[defi]{Corollary}
\newtheorem{rema}[defi]{Remark}
\newtheorem{cons}[defi]{Construction}
\newtheorem{ques}[defi]{Question}
\newtheorem{conj}[defi]{Conjecture}
\newtheorem{discussion}[defi]{Discussion}

\newcommand{\abs}{${\bar A}^*}
\newcommand{\qbs}{${\bar Q}^*}
\newcommand{\be}{\begin{enumerate}}
\newcommand{\ee}{\end{enumerate}}
\newcommand{\fany}{\rm for\ \ any\ \ }

\newcommand{\rb}{\overline{R}}
\newcommand{\mt}{\overline{M}}
\newcommand{\nt}{\overline{N}}
\newcommand{\nb}{\widetilde{N}}
\newcommand{\mb}{\widetilde{M}}
\newcommand{\m}{\bf {m}}


\def\cm{Cohen-Macaulay}
\def\wrt{with respect to\ }
\def\pni{\par\noindent}
\def\wma{we may assume without loss of generality that\ }
\def\Wma{We may assume without loss of generality that\ }
\def\ets{it suffices to show that\ }
\def\bwoc{by way of contradiction}
\def\iff{if and only if\ }
\def\st{such that\ }
\def\fg{finitely generated}


\def\a{\goth a}

\def\p{\mathbf p}

\def\isom{\thinspace \cong\thinspace}
\def\rtar{\rightarrow}
\def\rta{\rightarrow}
\def\l{\lambda}
\def\d{\Delta}

\def\alert#1{\smallskip{\hskip\parindent\vrule%
\vbox{\advance\hsize-2\parindent\hrule\smallskip\parindent.4\parindent%
\narrower\noindent#1\smallskip\hrule}\vrule\hfill}\smallskip}

\title{\bf Projective  equivalence of ideals in\\ Noetherian
integral domains  }

\author{William J. Heinzer, Louis J. Ratliff, Jr., and
David E. Rush}

\maketitle

\begin{abstract}
Let $I$ be a nonzero proper ideal in a Noetherian integral
domain $R$.
In this
paper we establish the existence of a
finite separable integral extension domain $A$ of $R$
and a positive integer $m$
such that all the Rees integers of $IA$ are equal to $m$.
Moreover, if $R$ has altitude one, then all the Rees integers
of $J$ = $\Rad(IA)$ are equal to one and
the ideals $J^m$ and $IA$
have the same integral closure. Thus
$\Rad(IA)$ =
$J$ is a projectively full radical ideal that is projectively
equivalent to $IA$.
In particular, if $R$ is Dedekind, then
there exists a Dedekind domain
$A$ having the
following properties: (i) $A$ is a finite separable integral
extension of $R$;  and,
(ii) there exists a radical ideal $J$ of $A$
and a positive integer  $m$  such that  $IA$ $=$ $J^m$.
In this case
the extension $A$ also has the property
that for each maximal
ideal $N$ of $A$ with $I \subseteq N$,
the canonical inclusion
$R/(N \cap R) \hookrightarrow A/N$ is an isomorphism,
and  the integer $m$
is a multiple of $[A_{(0)} : R_{(0)}]$.
\end{abstract}

\section{Introduction.}
All rings in this paper are commutative with a unit  $1$ $\ne$ $0$.
Let $I$ be a regular proper ideal of the Noetherian ring $R$, that
is, $I$ contains a regular  element of $R$ and $I \ne R$. An ideal
$J$ of $R$ is {\bf projectively equivalent} to $I$ if there exist
positive integers $m$ and $n$ such that $I^m$ and $J^n$ have the
same integral closure, that is, $(I^m)_a = (J^n)_a$, where  $K_a$
denotes the integral closure in  $R$ of an ideal $K$ of  $R$.  The
concept of projective equivalence of ideals and the study of ideals
projectively equivalent to $I$ was introduced by Samuel in \cite{S}
and further developed by Nagata in \cite{N1}.

Making use of interesting work of Rees in \cite{Re}, McAdam,
Ratliff, and Sally in \cite[Corollary 2.4]{MRS} prove  that the set
$\mathbf P(I)$ of integrally closed ideals projectively equivalent
to $I$ is discrete and linearly ordered with respect to  inclusion.
They also prove the existence of a fixed positive integer $d$ such
that for every ideal $J$ projectively equivalent to $I$,
$(J^d)_a = (I^n)_a$ for some positive integer $n$. If   $J$ and $K$
are  in $\mathbf P(I)$ and $m$ and $n$ positive integers, then
$(J^mK^n)_a \in \mathbf P(I)$. Thus  there is naturally associated
to $I$  a unique subsemigroup $S(I)$ of the additive semigroup of
nonnegative integers $\mathbb N_0$ such that $S(I)$ contains all
sufficiently large integers. A semigroup having these properties is
called a {\bf numerical semigroup}. The numerical semigroup $S(I)$
is an invariant of the projective equivalence class $\mathbf P(I)$
of $I$ in the sense that if $I$ is projectively equivalent to $J$,
then $S(I) = S(J)$, cf.
\cite[Remark~4.3]{CHRR}. It is observed in \cite[Remark 3.11]{CHRR2}
that every numerical semigroup  is realizable as $S(M)$ for an
appropriate local domain $(R,M)$.

The set $\mathbf P(I)$ is said to be {\bf projectively full} if $S(I) =
\mathbb N_0$, or equivalently,  if every element of $\mathbf P(I)$
is the integral closure of a power of the largest element $K$ of
$\mathbf P(I)$, i.e., every element of $\mathbf P(I)$ has the form
$(K^n)_a$, for some positive integer $n$. If this holds, then each
ideal $J$ in $R$ such that $J_a = K$ is said to be {\bf projectively
full}. A number of results about, and examples of, projectively full
ideals are given in \cite{CHRR}, \cite{CHRR2}, \cite{CHRR3}, and
\cite{CHRR4}.  Several characterizations of such ideals are given in
\cite[(4.11) and (4.12)]{CHRR}, and in \cite[Section 3]{CHRR2}
relations between projectively full ideals in  $R$  and in factor
rings of  $R$, localizations of  $R$, and extension rings of  $R$
are proved.

The set $\Rees I$ of Rees valuation rings of $I$ is a finite set of
rank one discrete valuation rings  (DVRs) that determine the
integral closure $(I^n)_a$ of $I^n$ for every positive integer $n$
and are the unique minimal set of DVRs having this property.
Consider the minimal primes $z$ of $R$ such that $IR/z$ is a proper nonzero
ideal. The set $\Rees I$ is the union of the sets $\Rees IR/z$. Thus
one is reduced to describing the set $\Rees I$ in the case where $I$
is a nonzero proper ideal of a Noetherian integral domain $R$.
Consider the Rees ring $\mathbf R = R[t^{-1}, It]$. The integral
closure $\mathbf R'$ of $\mathbf R$ is a Krull domain, so $W =
\mathbf R'_p$ is a DVR for each minimal prime $p$ of $t^{-1}\mathbf
R'$, and $V = W \cap F$, where $F$ is the field of fractions of $R$,
is also a DVR. The set $\Rees I$ of Rees valuation rings of $I$ is
the set of DVRs $V$ obtained in this way, cf.
\cite[Section~10.1]{SH}.

If  $(V_1, N_1), \ldots, (V_n, N_n)$ are the Rees valuation rings of
$I$, then the  integers $(e_1, \ldots, e_n)$, where $IV_i =
N_i^{e_i}$, are the {\bf Rees integers } of $I$.  Necessary and
sufficient conditions  for two regular proper ideals $I$ and $J$ to
be projectively equivalent are that (i) $\Rees I = \Rees J$ and (ii)
the Rees integers of $I$ and $J$ are proportional
\cite[Theorem~3.4]{CHRR}. If $I$ is integrally closed and each Rees
integer of $I$ is one, then $I$ is a projectively full radical
ideal.\footnote{There exist local domains $(R,M)$ for which $M$ is
not projectively full. A sufficient, but not necessary,  condition
in order that $I$ be projectively full is that the gcd of the Rees
integers of $I$ be one.  }

A main goal in the papers \cite{CHRR}, \cite{CHRR2}, \cite{CHRR3},
\cite{CHRR4}, and \cite{HRR},  and also in the present paper, is to
answer the following question:

\begin{ques}
\label{QUES} {\em Let  $I$  be a nonzero proper ideal in a
Noetherian domain  $R$.  Under what conditions does there exist   a
finite integral extension domain  $A$  of  $R$  such that  $\mathbf
P(IA)$ contains an ideal $J$  whose Rees integers are all equal to
one?

}
\end{ques}

Progress  is made on  Question \ref{QUES} in \cite{CHRR3}. To
describe this progress, let  $b_1,\dots,b_g$ be regular elements in
$R$ that generate  $I$ and for each positive integer $m$ $>$ $1$ let
$A_m$ $=$ $R[ {x_1},\dots, {x_g}]$ $=$ $R[ {X_1},\dots,
{X_g}]/({X_1}^m - b_1,\dots,{X_g}^m-b_g)$ and let $J_m$ $=$
$(x_1,\dots,x_g)A_m$. Then the main result in \cite{CHRR3}
establishes  the following:

\begin{theo}
\label{intro1}
Let  $R$  be a Noetherian ring, let  $I$  be a
regular proper ideal in  $R$, let  $b_1,\dots,b_g$ be regular
elements in  $I$  that generate  $I$, and let
$(V_1,N_1),\dots,(V_n,N_n)$ be  the Rees valuation rings of $I$.
Assume that:  {\rm (a)} $b_iV_j$ $=$ $IV_j$ ($=$ ${N_j}^{e_j}$, say)
for $i$ $=$ $1,\dots,g$  and $j$ $=$ $1,\dots,n$; and, {\rm (b)} the
greatest common divisor  $e$  of $e_1,\dots,e_n$ is a unit in  $R$.
Then $A_e$ $=$ $R[  {x_1},\dots, {x_g}]$ is a finite free integral
extension ring of $R$  and the ideal  $J_e$ $=$ $( {x_1},\dots,
{x_g})A_e$  is projectively full and projectively equivalent to
$IA_e$. Thus $\mathbf P(IA_e)$ $=$ $\mathbf P(J_e)$ is projectively
full.  Also, if $R$  is an integral domain and if $z$ is a
minimal prime ideal in  $A_e$, then
$((J_e+z)/z)_a$ is a projectively full ideal in  $A_e/z$
that is projectively equivalent to $(IA_e+z)/z$, so
$\mathbf P((IA_e+z)/z)$  is projectively full.
\end{theo}

We prove in
\cite[(3.19) and (3.20)]{HRR} that if either (i) $R$
contains an infinite field,  or (ii) $R$ is a local ring with an
infinite residue field, then it is possible to choose generators
$b_1, \ldots, b_g$ of $I$ that satisfy assumption (a) of
Theorem~\ref{intro1}. We prove in
\cite[(3.7)]{HRR} that if
``greatest common divisor'' is replaced with ``least common
multiple'', then the integral closure of the ideal   $J_e$  in
Theorem \ref{intro1} is a radical ideal with all Rees integers equal
to one.  Specifically:

\begin{theo}
\label{intro2}
With the notation of Theorem \ref{intro1}, assume
that: assumption {\rm (a)} of Theorem \ref{intro1} holds; and,
{\rm (b$^\prime$)}
the least common multiple  $c$  of $e_1,\dots,e_n$ is
a unit in  $R$. Then for each positive multiple $m$  of  $c$ that is
a unit in  $R$  the ideal $(J_m)_a$  is projectively full and
$(J_m)_a$  is a radical ideal that is projectively equivalent to
$IA_m$. Also, the Rees integers of  $J_m$ are all equal to one and
$x_iU$ is the maximal ideal of $U$  for each Rees valuation ring $U$
of $J_m$ and for $i$ $=$ $1,\dots,g$. Moreover, if  $R$ is an
integral domain and if $z$ is a minimal prime ideal in  $A_m$, then
$((J_m+z)/z)_a$ is a projectively full radical ideal that is
projectively equivalent to $(IA_m+z)/z$.
\end{theo}

Examples \cite[(3.22) and (3.23)]{HRR} show that condition
(b$^\prime$)
of Theorem
\ref{intro2}
is needed for the  proof of this result
given in
\cite{HRR}.
We show in \cite[(2.6)]{HRR} that every basis consisting
of regular elements of
$I$  can be used to find an integral extension ring
$A_m$  of  $R$  having a radical ideal  $J_m$  that
is projectively equivalent to  $IA_m$.  Specifically:

\begin{theo}
\label{intro3}
With notation as in Theorem \ref{intro1},
if $b_1,\dots,b_g$ are arbitrary
regular elements in $I$ that generate  $I$  and if $m$  is an integer
greater than or equal to  $\max(\{e_i \mid i = 1,\dots,n\})$,
then $J_m$  is projectively equivalent to
$IA_m$, $(J_m)_a$ $=$ $\Rad(J_m)$, and  $A_m/(J_m)_a$ $\cong$
$R/\Rad(I)$. Further, if  $R$  is an integral domain and if  $z$  is a
minimal prime ideal in  $A_m$, then $((J_m+z)/z)_a$ is a radical
ideal that is projectively equivalent to  $(IA_m+z)/z$.
\end{theo}

The main result in the present paper, Theorem \ref{prin}.2, answers
Question~\ref{QUES}  in the affirmative    for each nonzero proper
ideal $I$  in an arbitrary Noetherian integral domain  $R$  of
altitude one with no additional conditions; therefore the
conclusions of Theorems \ref{intro1} and \ref{intro2} are valid
without the assumption of conditions (a), (b), and (b$^\prime$) if
$R$  is a Noetherian integral domain of altitude one.
In particular, Theorem
\ref{prin}.2
shows that these conclusions
hold for the
above mentioned
examples
\cite[(3.22) and (3.23)]{HRR}.

A classical theorem of Krull is an important tool in our work. By
successively applying this theorem of Krull,  we construct a finite
integral extension domain $A$ of $R$ such that $H$ $=$ $\Rad(IA)$ is
a projectively full radical ideal
 that is projectively equivalent to  $IA$.  Moreover,
the Rees integers of $H$ are all equal to one. If, in addition, $R$
is integrally closed, then $A$ is the integral closure of $R$ in a
finite separable algebraic field extension and $H^m = IA$, where $m$
is a multiple of $[A_{(0)} : R_{(0)}]$; and for each maximal ideal
$N$ of $A$ with $I \subset N$, the canonical inclusion map $R/(N
\cap R) \hookrightarrow A/N$ is an isomorphism.

In Section~3
we consider the question of extending Theorem~\ref{prin}.2 to the case
of regular principal ideals  $bR$  of a Noetherian domain  $R$  of altitude
greater than one. A complicating factor here is the possibility that
$\Rad(bA)$  may have embedded asymptotic prime divisors.  In Section~4 we
present an
application that partially extends Theorem~\ref{prin}.1 to certain
finite sets of ideals.

Our notation is mainly as in Nagata \cite{N2}, so, for example,
the term {\bf altitude}  refers to what is often called  dimension or Krull
dimension, and a {\bf basis} for an ideal is a set of elements that generate
the ideal.

\section{Finite integral extensions of a Noetherian domain.}

To prove our main result, we use a theorem of Krull; before stating
Krull's Theorem, we recall the following terminology from
\cite{Gilmer}.

\begin{defi}
\label{consist} {\em Let $(V_1,N_1), \ldots, (V_n,N_n)$ be distinct
rank one discrete valuation domains of a field $F$ and for  $i$ $=$
$1,\dots,n$  let $K_i$ $=$ $V_i/N_i$ denote the residue field of
$V_i$. Let $m$ be a positive integer. By an {\bf m-consistent system
for} $\{ V_1, \ldots, V_n \}$, we mean a collection of sets $S$ =
$\{ S_1, \ldots, S_n \}$ satisfying the following conditions:

(1) $S_i$ = $\{ (K_{i,j}, f_{i,j},e_{i,j}) \mid j = 1, \dots, s_i \}$, where
$K_{i,j}$ is a simple algebraic field extension of $K_i$ and  $s_i,
f_{i,j},e_{i,j} \in \mathbb N_+$  (the set of positive integers).

(2) For each $i$, the sum $\sum_{j=1}^{s_i} e_{i,j}f_{i,j}$
= $m$. }
\end{defi}

\begin{defi}
\label{real} {\em The  $m$-consistent system $S$ as in
Definition~\ref{consist} is said to be {\bf realizable} if  there
exists a separable algebraic extension field $L$ of $F$ such that:

(a) $[L : F]$ = $m$.

(b) For $1 \leq i \leq n$, $V_i$ has exactly $s_i$ extensions
$V_{i,1}, \ldots, V_{i,s_i}$ to $L$.

(c) The residue field of $V_{i,j}$ is $K_i$-isomorphic to $K_{i,j}$,
the residue field extension  $K_{i,j}$  of  $K_i$  has degree  $f_{i,j}$
(so $[K_{i,j} : K_i]$ $=$ $f_{i,j}$),
and the ramification index of $V_{i,j}$ over $V_i$ is $e_{i,j}$  (so
$N_iV_{i,j}$ $=$ ${N_{i,j}}^{e_{i,j}}$).

\noindent If  $S$  and  $L$  are as above, we say the field  $L$
{\bf{realizes}} $S$  or that  $L$  is a {\bf{realization}} of  $S$.
}
\end{defi}

\begin{theo}
\label{GK}
{\em (Krull \cite{WKrull}):}
Let $(V_1,N_1), \ldots, (V_n,N_n)$
be distinct rank one discrete
valuation domains of a field $F$
with $K_i$ $=$ $V_i/N_i$
for  $i$ $=$
$1,\dots,n$, let  $m$ be a positive integer, and let
$S$ = $\{ S_1, \ldots, S_n \}$ be an $m$-consistent system for
$\{ V_1, \ldots, V_n \}$
with
$S_i$ = $\{ (K_{i,j}, f_{i,j},e_{i,j}) \mid j = 1, \dots, s_i \}$
for  $i$ $=$ $1,\dots,n$.
Then $S$ is realizable if one of the following
conditions is satisfied.

(i) $s_i$ = $1$ for at least one $i$.

(ii) $F$ has at least one rank one discrete
valuation domain  $V$ distinct from
$V_1, \ldots, V_n $.

(iii) For each monic polynomial $X^t + a_1X^{t-1} + \cdots + a_t$
with $a_i \in \cap_{i=1}^n V_i$ = $D$, and for each $h \in \mathbb
N$, there exists an irreducible separable polynomial $X^t +
b_1X^{t-1} + \cdots + b_t \in D[X]$ with $b_l - a_l \in {N_i}^h$ for
each $l$ = $1, \ldots, t$ and  $i$ = $1, \ldots, n$.
\end{theo}

Observe that condition (i) of Theorem \ref{GK} is a property of the
$m$-consistent system $S$ = $\{ S_1, \ldots, S_n \}$, whereas
condition (ii) is a property of the family of rank one discrete
valuation domains of the field $F$,  and condition (iii) is a
property of  the family $(V_1,N_1), \ldots, (V_n,N_n)$.

\begin{rema}
\label{Dedekind}
{\em
Let $D$ be a Dedekind domain
with quotient field $F \neq D$,
let
$M_1,\dots,M_n$  be distinct maximal ideals of $D$,
let
$I$ = ${M_1}^{e_1} \cdots {M_n}^{e_n}$
be an ideal in $D$,
where $e_1,\dots,e_n$
are positive integers,
and
let $S$ = $\{ S_1, \ldots, S_n \}$  be a
realizable $m$-consistent system  for
$\{ D_{M_1},\dots,D_{M_n} \}$, where
$S_i$ = $\{ (K_{i,j},f_{i,j}, e_{i,j}) \mid j = 1, \dots, s_i \}$
for  $i$ $=$ $1,\dots,n$.
Let  $L$  be a field that realizes  $S$  and let  $E$  be
the integral closure of  $D$  in  $L$.  Then:

\noindent
{\bf{(\ref{Dedekind}.1)}}
$L$  has rank one
discrete valuation domains
$(V_{i,1}, N_{i,1}), \dots,(V_{i,s_i},N_{i,s_i})$
such that
for each  $i,j$:
$V_{i,j} \cap F$ $=$ $D_{M_i}$;
$V_{i,j}/N_{i,j}$ is
$D/M_i$-isomorphic to $K_{i,j}$;
$[K_{i,j} : (D/M_i)]$ $=$ $f_{i,j}$; and,
$M_iV_{i,j}$ $=$ ${N_{i,j}}^{e_{i,j}}$.
Also, for  $i$ $=$ $1,\dots,n$,
$V_{i,1},  \ldots, V_{i,s_i}$
are all of the extensions of  $D_{M_i}$
to  $L$.

\noindent
{\bf{(\ref{Dedekind}.2)}}
$E$
is a Dedekind domain
that is a finite separable integral extension domain of  $D$, and
$IE$ = ${M_1}^{e_1} \cdots {M_n}^{e_n}E$
= ${P_{1,1}}^{e_1e_{1,1}} \cdots {P_{n,s_n}}^{e_ne_{n,s_n}}$,
where  $P_{i,j}$ $=$ $N_{i,j} \cap E$
for $i$ $=$ $1,\dots,n$ and
$j$ = $1, \ldots, s_i$.
}
\end{rema}

\begin{proof}
(\ref{Dedekind}.1) follows immediately from (a) - (c) of
Definition~\ref{real}.

For (\ref{Dedekind}.2), $E$ is a Dedekind domain, by \cite[Theorem
19, p. 281]{ZS1}, and  $E$  is a finite separable integral extension domain
of
$D$, by \cite[Corollary 1, p. 265]{ZS1}, since  $L$  is a finite
separable algebraic extension field of  $F$.  Also, $V_{i,j}$ $=$
$E_{P_{i,j}}$, so $IV_{i,j}$ $=$ $(IE)V_{i,j}$ $=$
$(ID_{M_i})V_{i,j}$ $=$ $({M_i}^{e_i}D_{M_i})V_{i,j}$ $=$
$({M_i}{V_{i,j}})^{e_i}$ $=$ ${N_{i,j}}^{e_ie_{i,j}}$. Since the
ideals  $P_{i,j}$  are the only prime ideals in  $E$  that lie over
$M_i$  (for  $i$ $=$ $1,\dots,n$ and $j$ = $1, \ldots, s_i$) and
since the  $P_{i,j}$  are comaximal, it follows that $IE$ $=$
${P_{1,1}}^{e_1e_{1,1}} \cdots {P_{n,s_n}}^{e_ne_{n,s_n}}$.
\end{proof}

We use the following two lemmas in the proof of Theorem \ref{prin}.

\begin{lemm}
\label{inductnew} Let  $D$ be a Dedekind domain  and let  $I$ =
${M_1}^{e_1} \cdots {M_n}^{e_n}$  ($n$ $>$ $1$) be an irredundant
primary decomposition of  a nonzero proper ideal $I$  in  $D$.
Assume that  the  integers $e_i$ have no common factor $d > 1$. Let
$p$ be a prime integer dividing at least one of the $e_i$.
 Then there
exists a Dedekind domain  $E_1$  that is a finite separable integral
extension domain  of  $D$  with an ideal  $J_1$  such that:
${J_1}^{p^h}$ $=$ $IE_1$  for some positive integer  $h$;  and, if
${J_1}$ $=$ ${N_1}^{c_1} \cdots {N_g}^{c_g}$ is an irredundant
primary decomposition of  ${J_1}$, then $\prod_{j=1}^g c_j$ has
fewer distinct prime integer factors  than  does $\prod_{j=1}^n
e_j$. Moreover, the canonical inclusion map $D/(N_i \cap D)
\hookrightarrow E_1/N_i$ is an isomorphism for each $i \in \{1,
\ldots, g\}$.

\end{lemm}

\begin{proof}
Let $e_i$ = $p^{h_i}d_i$, where  $p \not| d_i$, and $h_i \ge 0$.  We
may assume that the $e_i$ are ordered so that  $h_1 \geq h_2 \geq
\cdots \geq h_n$. Our hypotheses imply that $h_1 > 0$ and $h_n = 0$.
Let $S$ $=$ $\{S_1,\dots,S_n\}$ with $S_i$ = $\{(K_{i,j},1,p^{h_1-h_i})
\mid j = 1,\dots, p^{ h_i}\}$.
We show that  $S$  is a realizable
$p^{h_1}$-consistent system for  $\{D_{M_1},\dots,D_{M_n}\}$.
Observe that $\sum_{j=1}^{s_i} e_{i,j}f_{i,j}$
$=$ $\sum_{j=1}^{p^{h_i}} p^{h_1-h_i} \cdot 1$
 $=$ $ p^{h_1}$. Therefore $S$  is a
$p^{h_1}$-consistent system. Since  $s_n$ $=$ $p^{h_n}$ = $p^0$ =
$1$, $S$ is realizable, by Theorem \ref{GK}(i). Therefore, by Remark
\ref{Dedekind} (especially (\ref{Dedekind}.2)), the integral closure
$E_1$  of  $D$  in a realization $L$  of  $S$ for
$\{D_{M_1},\dots,D_{M_n}\}$  is a Dedekind domain such that $IE_1$ =
$\prod_{i=1}^n ({M_i}^{e_i}E_1)$ =
$$\prod_{i=1}^n (
{N_{i,1}}^{e_ie_{i,1}}
\cdots
{N_{i,s_i}}^{e_ie_{i,s_i}})
=
\prod_{i=1}^n (
{N_{i,1}}^{(p^{h_i}d_i)(p^{h_1-h_i})}
\cdots
{N_{i,s_i}}^{(p^{h_i}d_i)(p^{h_1-h_i})}
) =
{J_1}^{p^{h_1}},
$$ where
${J_1}$ = $\prod_{i=1}^n ( {N_{i,1}}^{d_i} \cdots {N_{i,s_i}}^{d_i}
)$. Also, $\prod_{i=1}^n {d_i}^{s_i}$ $=$ $\prod_{i=1}^n
{d_i}^{p^{h_i}}$ has fewer distinct prime integer factors  than does
$\prod_{i=1}^n e_i$.
Finally, since all  $f_{i,j}$  are equal to one, it follows
that $E_1/N$  $\cong$ $D/(N \cap D)$  for all maximal
ideals  $N$  of  $E_1$ that contain  $I$.
\end{proof}

\begin{lemm}
\label{prin.reduction.lemma} Let $R$ be a  Noetherian domain of
altitude one with quotient field $F$ and let $I$  be a nonzero
proper ideal in $R$. Let $L$ be a finite algebraic extension field
of $F$, and let $E$ denote the integral closure of $R$ in $L$. If
there exist distinct maximal ideals $N_1, \ldots, N_n$ of $E$ and
positive integers  $k_1,\dots,k_n,h$ such that $IE$ = $({N_1}^{k_1}
\cdots {N_n}^{k_n})^h$, then there exists a finite integral
extension domain $A$ of $R$ with quotient field $L$ and an ideal $H$
of $A$
 such that
 \begin{enumerate}
 \item[{\rm(i)}]
  $H$ has Rees valuation rings  $E_{N_{1}},\dots,E_{N_n}$ with
corresponding Rees integers  $k_1, \ldots, k_n$; and,

\item[{\rm (ii)}]  $(H^h)_a$
= $(IA)_a$.
\end{enumerate}
\end{lemm}

\begin{proof}
If  $E$  is a finite  $R$-module, then let $A$ $=$ $E$  and  $H$ $=$
${N_1}^{k_1} \cdots {N_n}^{k_n}$. Otherwise, let  $A_0 $  be a
subring of $E$ that is a finite integral extension domain of  $R$
and that has quotient field $L$. For $i$ $=$ $1,\dots,n$ let $G_i
\subseteq E$ be a finite set such that $G_iE$ = $N_i$, let $A$ =
$A_0[G_1 \cup \cdots \cup  G_n]$, and for $i$ $=$ $1,\dots,n$ let
$P_{i}$ $=$ $N_{i} \cap A$.  Then  $A$  is a finite integral
extension domain of  $R$,  $E$  is the integral closure of  $A$  in
its quotient field $L$, and  $P_{i}$  is a maximal ideal in  $A$
such that  $P_{i}E$ $=$ $N_{i}$ for $i$ $=$ $1,\dots,n$. Let $H$ =
${P_1}^{k_1} \cdots {P_n}^{k_n}$. Since
$$
HE \cap A = [({P_1}^{k_1} \cdots {P_n}^{k_n})E] \cap A =
[({P_1}^{k_1}E) \cap \cdots \cap ({P_n}^{k_n}E)] \cap A =
[{N_1}^{k_1} \cap \cdots \cap {N_n}^{k_n}] \cap A, $$ our hypotheses
imply that $H$ has Rees valuation rings $E_{N_{1}},\dots,E_{N_n}$
with corresponding Rees integers  $k_1, \ldots, k_n$. Also $(H^h)_a$
= $(H^hE)_a \cap A$ = $(H^hE) \cap A$ = $({N_1}^{k_1h} \cdots
{N_n}^{k_nh}) \cap A$ = $IE \cap A$ = $(IE)_a \cap A$ = $(IA)_a $.
\end{proof}

We also use the following  well-known fact concerning the Rees
valuation rings  of an ideal (cf.  the proof of \cite[Theorem
2.5]{CHRR3}).

\begin{rema}
\label{fact} {\em Let  $I$  be a nonzero proper ideal in a
Noetherian  integral domain  $R$, let  $\Rees I$  $=$
$\{V_1,\dots,V_n\}$, and let $A$  be a finite integral extension
domain of  $R$.  Then $\Rees IA$ $=$ $\{V_{1,1}, \dots,V_{n,c_n}\}$,
where for each $i \in \{1, \ldots, n\}$,  $V_{i,1},\dots,V_{i,c_i}$
are all the extensions of $V_i$ to the quotient field of  $A$.

}
\end{rema}

Theorem \ref{prin}.2, answers Question~\ref{QUES}  in the
affirmative for each nonzero proper ideal  in an arbitrary
Noetherian integral domain   of altitude one with no additional
conditions.

\begin{theo}
\label{prin} Let $I$ be a nonzero proper ideal in a Noetherian
integral domain $R$.
\begin{enumerate}
\item  There exists a finite separable integral extension domain $A$
of $R$ and a positive integer $m$ such that all the Rees integers of
$IA$ are equal to $m$.
\item If $R$ has altitude one, then
there exists a finite integral extension domain  $A$  of  $R$  such
that $\mathbf P(IA)$ contains an ideal $H$  whose Rees integers are
all equal to one. Therefore $H$ $=$ $\Rad(IA)$  is a projectively
full radical ideal that is projectively equivalent to  $IA$.
\end{enumerate}
\end{theo}

\begin{proof}
For part 2,
if $R$ has altitude one, then the integral closure of $R$ is a
Dedekind domain $D$ and there exist distinct maximal ideals
$M_1,\dots,M_n$ ($n$ $\ge$ $1$) in $D$  and positive integers
$e_1,\dots,e_n$  such that $ID$ = ${M_1}^{e_1} \cdots {M_n}^{e_n}$.
The $D_{M_i}$ are the Rees valuation rings of $I$ and
$e_1,\dots,e_n$ are the Rees integers of $I$. If either $n$ $=$ $1$
or $e_1$ $=$ $\cdots$ $=$ $e_n$, then the conclusions of part 2 follow from
Lemma \ref{prin.reduction.lemma} with $L$ = $F$, so we may assume
that $n$ $>$ $1$ and that not all the $e_i$ are equal. Let $d$  be
the greatest common divisor of  $e_1,\dots, e_n$. Then the ideal
$I_0$ $=$ ${M_1}^{\frac{e_1}{d}} \cdots {M_n}^{\frac{e_n}{d}}$  is
such that  $ID$ $=$ ${I_0}^d$, so the ideal  $I_0$  may be used in
place of  $ID$.  Thus we may assume that the $e_i$ have no common
factor $d > 1$. Let $k$ be the number of distinct prime integers
dividing $\prod_{j=1}^n e_j$. By induction on $k$,  it suffices to
show that there exists a finite integral extension domain $A$ of
$R$, an ideal $H$ of $A$, and a positive integer $h$ such that:
$(H^h)_a$ = $(IA)_a$; $H$ has Rees integers $c_1, \ldots, c_g$; and,
there are at most $k - 1$ distinct prime integers  dividing
$\prod_{j=1}^g c_j$. Therefore Theorem~\ref{prin}.2  follows from
Lemmas \ref{inductnew}
and  \ref{prin.reduction.lemma}.

For the proof of part 1, let $\Rees I = \{(V_1,N_1), \ldots,
(V_n,N_n)\}$, let $e_1, \ldots, e_n$ be the Rees integers of $I$,
 and let $D = V_1 \cap \cdots \cap V_n$. Then $D$ is a Dedekind
 domain with maximal ideals $M_1, \ldots, M_n$ and $D_{M_i} = V_i$
 for each $i$ with $1 \le i \le n$. Also $ID = M_1^{e_1} \cdots
 M_n^{e_n}$. If either $n = 1$ or $e_1 = \cdots = e_n$, then the
 assertion of part 1 is obvious. Thus we may assume that $n > 1$ and
 that not all the $e_i$ are equal. The argument in the paragraph
 above for part 2 implies that there exists a finite separable
 algebraic field extension $L$ of the quotient field $F$ of $D$ such
 that if $E$ is the integral closure of $D$ in $L$, then $IE = J^m$,
 where $J$ is a radical ideal of $E$. There exists $\theta \in L$
 such that $L = F[\theta]$ and there exists a nonzero $r \in R$ such
 that $r\theta$ is integral over $R$. Let $A = R[r\theta]$.
 Remark~\ref{fact} implies that each of the Rees integers of $IA$ is
 $m$.
\end{proof}

\begin{coro} \label{2.8}
Let $I$ be a nonzero proper ideal in a Dedekind domain $D$. There
exists a Dedekind domain $E$ having the following properties: (i)
$E$ is a finite separable integral extension of $D$; and, (ii) there
exists a radical ideal $J$ of $E$ and a positive integer  $m$ such
that $IE$ $=$ $J^m$. Therefore $J$ is a projectively full radical
ideal that is projectively equivalent to $IE$, and the Rees integers
of $J$ are all equal to one. The extension $E$ also has the property
that for each maximal ideal $N$ of $E$ with $I \subseteq N$, the
canonical inclusion $D/(N \cap D) \hookrightarrow E/N$ is an
isomorphism, and  $m$ is a multiple of $[E_{(0)} : D_{(0)}]$.
\end{coro}

\begin{proof} Everything but the last sentence of
Corollary~\ref{2.8} is immediate from Theorem~\ref{prin}.2. The
application of Lemma~\ref{inductnew} to the integral closure $D$ of
the Noetherian domain $R$ of Theorem~\ref{prin}.2 implies that $D/(N
\cap D) \hookrightarrow E/N$ is an isomorphism. That $m$ is a
multiple of $[E_{(0)} : D_{(0)}]$ follows from Remark~\ref{misc}.3.
\end{proof}

In Lemma~\ref{induct}, we give a different consistent system for
$D_{M_1},\dots,D_{M_n}$ that may  be used in place of Lemma
\ref{inductnew} to inductively complete an alternative  proof of
Theorem \ref{prin}. The  proof of Lemma~\ref{induct} is described
more fully in Remark \ref{misc}.1.  (Concerning the hypothesis
``$k$ $>$ $1$'' in Lemma \ref{induct}, if  $k$ $=$ $0$, then
$I$ $=$ $M_1 \cdots M_n$  is a raddical ideal and the lemma holds
with  $E_1$ $=$ $D$, $J_1$ $=$ $I$, and  $h$ $=$ $1$.)

\begin{lemm}
\label{induct}
Let  $D$ be a Dedekind domain, let  $I$ =
${M_1}^{e_1} \cdots {M_n}^{e_n}$  ($n$ $>$ $1$) be an irredundant
primary decomposition of  a nonzero proper ideal $I$  in  $D$, and
assume that $e_i$ $>$ $1$ for at most  $k$ $(1 \leq k \leq n)$
of the
integers  $e_i$.
Then there exists a Dedekind domain  ${E_1}$  that
is a finite separable integral extension domain  of  $D$  with an
ideal  ${J_1}$  such that:  ${J_1}^h$ $=$ $I{E_1}$  for some
positive integer  $h$;  and, if  ${J_1}$ $=$ ${N_1}^{c_1} \cdots
{N_g}^{c_g}$ is an irredundant primary decomposition of  ${J_1}$,
then $c_j$ $>$ $1$  for at most  $k-1$  of the integers
$c_1,\dots,c_g$. Moreover, the canonical inclusion map $D/(N_i \cap
D) \hookrightarrow E_1/N_i$ is an isomorphism for each $i \in \{1,
\ldots, g\}$.

\end{lemm}

\begin{proof}
Let $S$ $=$ $\{S_1,\dots,S_n\}$,
where:   $S_1$ $=$ $\{(K_{1,j},1,1) \mid j = 1,\dots,e_1\}$;
and, for $i$ $=$ $2,\dots,n$,
$S_i$ $=$ $\{(K_{i,1},1,e_1)\}$.
Then a proof similar to the proof of
Lemma \ref{inductnew} shows that
$S$  is
a realizable $e_1$-consistent
system for  $\{D_{M_1},\dots,D_{M_n}\}$
and that  $I{E_1}$ $=$
${{J_1}}^{e_1}$, where
${J_1}$ $=$ $({N_{1,1}} \cdots {N_{1,e_1}}) {N_{2,1}}^{e_2} {N_{3,1}}^{e_3}
\cdots {N_{n,1}}^{e_n}$.
Finally, since all  $f_{i,j}$  are equal to one, it follows
that $E_1/N$  $\cong$ $D/(N \cap D)$  for all maximal
ideals  $N$  of  $E_1$ that contain  $I$.
\end{proof}

\begin{rema}
\label{misc}
{\bf{(\ref{misc}.1)}} {\em In Theorem \ref{prin}.2,
assume that the exponents $e_1,\dots,e_n$  are arranged so that
$e_i$ $>$ $1$  if and only if  $i$ $\in$ $\{1,\dots,k\}$, where  $k$
$\le$ $n$. If we successively carry out the separate steps of the
induction in the proof of Theorem \ref{prin}.2 using Lemma \ref{induct},
then we get a chain of rings  $D$ $=$ $E_0$ $\subset$ $E_1$
$\subset$ $\cdots$ $\subset$ $E_{k}$ $=$ $E$, where each  $E_{i}$
($i$ $=$ $1,\dots,k$)  is a Dedekind domain that is a finite
separable integral extension of $E_{i-1}$  and for which $(E_i)_U$
has exactly $e_i -1$  more maximal ideals than  $(E_{i-1})_U$, where
$U$ $=$ $D \setminus (M_1 \cup \cdots \cup M_n) $. In fact, for $i$
$=$ $1,\dots,k$, $E_{i}$  is obtained as the integral closure of
$E_{i-1}$  in a realization $L_{i}$  of a realizable
$e_{i}$-consistent system   $S^{(i)}$  for $\{(E_{i-1})_N \mid N \in
\mathbf N(E_{i-1})\}$, where  $\mathbf N(E_{i-1})$ $=$ $\{N \mid N$
is a maximal ideal in  $E_{i-1}$  and $N \cap D$ $\in$ $\{M_1
,\dots,M_n\}$. Here, the $e_{i}$-consistent system $S^{(i)}$
completely splits (into  $e_{i}$ components) the unique maximal
ideal in  $E_{i-1}$ that contracts in  $D$  to $M_{i}$, and it
completely ramifies (of index  $e_{i}$) all the remaining maximal
ideals in  $\mathbf N(E_{i-1})$, so $E_i/N$ $\cong$ $E_{i-1}/(N \cap
E_{i-1})$ $\cong$ $D/(N \cap D)$  for all $N$  in $\mathbf N(E_i)$.
($S^{(i)}$  is realizable, by Theorem \ref{GK}(i), since for all but
one  $N$ $\in$ $\mathbf N(E_{i-1})$, the corresponding  component
${S^{(i)}}_j$  of  $S^{(i)}$ contains a single ordered triple
$(E_{i-1}/N,1,e_i)$.) Therefore: (a)  exactly $e_{i}$  of the maximal
ideals in  $E_i$ contract in  $E_{i-1}$  to the unique maximal ideal
in  $E_{i-1}$  that contracts in  $D$  to $M_{i}$; and, (b)  the
remaining maximal ideals in  $\mathbf N(E_i)$ are in one-to-one
correspondence with the remaining $e_1+ \cdots +e_{i-1}+(n-i)$
maximal ideals in  $\mathbf N(E_{i-1})$. Further, for each  maximal
ideal  $N$  of  (b) it holds that  $E_i/N$ $\cong$ $D/(N \cap D)$ and
$(N \cap E_{i-1})(E_i)_N$ $=$ $N^{e_{i}}(E_i)_N$, while for the
$e_{i}$ maximal ideals $N$  of (a) it holds that  $E_i/N$ $\cong$ $D/(N \cap
D)$ and  $(N \cap E_{i-1})(E_i)_N$ $=$ $N(E_i)_N$. It follows that,
in  $E$ $=$ $E_k$, $\mathbf N(E)$ has exactly $e_1 + \cdots + e_k +
(n-k)$  maximal ideals, and of these, exactly $e_i$ of them contract
in  $D$ to  $M_i$  for  $i$ $=$ $1,\dots,n$.  Also, if  $N$ is a
maximal ideal in  $E$  and  $N \cap D$ $=$ $M_i$  (with $i$ $\in$
$\{1,\dots,n\}$), and if  ${e_i}^*$ $=$ $\frac{e_1 \cdots
e_n}{e_i}$, then  $N$ $\in$ $\mathbf N(E_k)$ ($=$ $\mathbf N(E)$),
$ME_n$ $=$ $N^{{e_i}^*}E_N$, and  $E/N$ $\cong$ $D/M_i$. It therefore
follows that:
\newline
 {\it{{\rm (*1)}  the quotient field  $L$ $=$ $L_k$ of  $E$  is a
realization of the realizable $e_1 \cdots e_n$-consistent system  $S
$ $=$ $\{S_1,\dots,S_n\}$ for  $\{D_{M_1},\dots,D_{M_n}\}$, where:
for $i$ $=$ $1,\dots,k$, $S_i$ $=$ $\{(K_{i,j},1,{e_i}^*) \mid
j = 1,\dots,e_i\}$; and, for  $i$
$=$ $k+1,\dots,n$, $S_i$ $=$ $\{(K_{i,1},1,e_1 \cdots e_n) \}$,\footnote{If
at least one of the integers
$e_1,\dots,e_n$  is one, then it follows from Theorem \ref{GK}(i)
that  $S$  is a realizable $e_1 \cdots e_n$-consistent system, so
(*1) readily follows.  However, if $e_i$ $\ne$ $1$  for  $i$ $=$
$1,\dots,n$, then it is only by this ``composition'' of realizable
consistent systems
that we are able to show that  $S$  is realizable, and thereby find
a finite integral extension domain $E$  of  $D$  for which  $IE$  is
the power of a radical ideal of  $E$.  This idea of composing
realizable consistent systems is further developed in \cite{HRR2}.}
so  $IE$ $=$ $J^{e_1 \cdots e_n}$,
where  $J$ $=$ $\cap \{N \mid N \in \mathbf N(E)\}$}}
(since  $IE_N$ $=$ $N^{e_1 \cdots e_n}E_N$
for each maximal ideal  $N$  in  $\mathbf N(E)$),
hence  $JE_U$  is the Jacobson radical of  $E_U$.

\noindent
{\bf{(\ref{misc}.2)}}
Assume that  $I$ $=$ ${M_1}^{e_1} \cdots {M_n}^{e_n}$,
that no prime integer divides each  $e_i$, and that
the least common multiple of  $e_1,\dots,e_n$  is
$d$ $=$ ${p_1}^{m_1} \cdots {p_k}^{m_k}$, where
$p_1,\dots,p_k$  are distinct prime integers
and  $m_1,\dots,m_k$  are positive integers.
Then it follows as in (\ref{misc}.1) that
if we successively carry out the separate steps of the induction in
the proof of
Theorem \ref{prin}.2
using Lemma \ref{inductnew}, then
we get a chain of rings  $D$ $=$ $E_0$ $\subset$ $E_1$ $\subset$
$\cdots$ $\subset$ $E_{k}$ $=$ $E$, where each  $E_{i}$  ($i$
$=$ $1,\dots,k$)  is a
Dedekind domain that is obtained as the
integral closure of  $E_{i-1}$  in a realization
$L_{i}$  of a realizable  ${p_i}^{m_i}$-consistent
system   $S^{(i)}$  for  $\{(E_{i-1})_N \mid N
\in \mathbf N(E_{i-1})\}$, where  $\mathbf N(E_{i-1})$
$=$ $\{N \mid N$  is
a maximal ideal in  $E_{i-1}$  and  $N \cap D$
$\in$ $\{M_1 ,\dots,M_n\}\}$.
It therefore follows that:
\newline
{\it{{\rm (*2)}  the quotient field  $L$ $=$ $L_k$
of  $E$  is a realization of the realizable $d$-consistent
system  $S $ $=$ $\{S_1,\dots,S_n\}$ for  $\{D_{M_1},\dots,D_{M_n}\}$,
where: for  $i$ $=$ $1,\dots,n$,
$S_i$ $=$ $\{(K_{i,j},1,\frac{d}{e_i}) \mid j = 1,\dots,e_i\}$, so
$IE$ $=$ $J^d$, where  $J$ $=$
$\prod_{i=1}^n (M_{i,1} \cdots M_{i,e_i})$}  with}  $\mathbf N(E)$
$=$ $\{N_{1,1},\dots,N_{n,e_n}\}$

\noindent {\bf{(\ref{misc}.3)}} It follows from the last part of
(\ref{misc}.1) that, in Corollary~\ref{2.8}, the extension domain
$E$ of $D$  and the integer  $h$ such that  $IE$ $=$ $J^h$ can be
chosen such that:  $h$ $=$ $e_1 \cdots e_n$; and, the quotient field
$L$ of  $E$  is a realization of an  $h$-consistent system for the
Rees valuation rings of  $I$. And it follows from the last part of
(\ref{misc}.2) that, in Corollary~\ref{2.8}, if  $d$  is the
greatest common divisor of  $e_1,\dots,e_n$, if  $I_0$ $=$
${M_1}^{\frac{e_1}{d}} \cdots {M_n}^{\frac{e_n}{d}}$ (so  $I_0$  is
projectively equivalent to  $I$), and if  $c$  is the least common
multiple of $\frac{e_1}{d}, \dots, \frac{e_n}{d}$, then the
extension domain $E$  of  $D$  and the integer  $h$ such that $I_0E$
$=$ $J^h$  can be chosen such that:  $h$ $=$ $c$; and, the quotient
field  $L$ of  $E$  is a realization of an $h$-consistent system for
the Rees valuation rings of  $I$  (and of $I_0$).

}
\end{rema}

\section{Principal ideals and projective equivalence in \\ finite integral
extensions.}

In this section we consider the question of an extension of Theorem
\ref{prin}.2 to regular principal ideals of a Noetherian integral
domain of altitude greater than one.

\begin{discussion} \label{3.1} {\em
Let  $b$  be a nonzero  nonunit in a Noetherian integral domain $R$,
let $R'$ be the integral closure of  $R$  in its quotient field $F$,
let $p_1,\dots,p_n$  be the height-one  prime ideals in  $R'$ that
contain  $bR'$, and let ${{p_1}}^{(e_1)} \cap \cdots \cap
{{p_n}}^{(e_n)}$ (symbolic powers) be an irredundant  primary
decomposition of  $bR'$. It follows (see, for example,
\cite[(2.3)]{CHRR2}) that the rings  $V_i$ $=$ ${R'}_{{p_i}}$  ($i$
$=$ $1,\dots,n$)  are the Rees valuation rings of  $bR$. Let  $D$
$=$ ${R'}_U$, where  $U$ $=$ $R' \setminus  ({p_1} \cup \cdots \cup
{p_n})$. Theorem \ref{prin}.2 implies that there exists a finite
separable algebraic extension field $L$ $=$ $F[\theta]$ of  $F$ such
that the integral closure  $E$  of  $D$  in  $L$ is a Dedekind
domain having  a radical ideal  $J$ such that $bE$ $=$ $J^m$ for
some positive integer  $m$. If  $\altitude(R)$ $=$ $1$, then $J \cap
R'[\theta]$ is a radical ideal that is projectively full and
projectively equivalent to $bR'[\theta]$, by Theorem \ref{prin}.2 and
its proof. Thus  it seems at least plausible that this may also hold
when $\altitude(R)$ $>$ $1$. However, a complication in higher
altitude is that  powers of $J \cap R'[\theta]$ may have embedded
asymptotic prime divisors, as the following example shows.
}
\end{discussion}

\begin{exam}
\label{failure} {\em Let $k$  be a field, let  $x,y$  be independent
indeterminates, let $R$ $=$ $k[[x^2, xy, y^2]]$, and let $P$ $=$
$(x^2, xy)R$. Then  $R$  is an integrally closed local domain of
altitude two, the regular local ring  $A$ $=$ $k[[x,y]]$ is a finite
integral extension domain of  $R$, and $P$ $=$ $xA \cap R$ is the
radical of the principal ideal $bR$ $=$ $x^2R$, so  $V$ $=$
${R}_{P}$ is the only Rees valuation ring of  $bR$. Also, $N$ $=$
$PV$ $=$ $(x^2,xy)V$  and $\frac{x}{y}$ $=$ $(xy) \cdot
\frac{1}{y^2}$ $\in$ $N$, so  $x^2$ $=$ $(xy) \cdot \frac{x}{y}$
$\in$ $N^2$, so it follows that  $N$ $=$ $xyV$ $=$ $\frac{x}{y}V$,
hence  $N^2$ $=$ $bV$, so  $N$  is a radical ideal that is
projectively equivalent to  $bV$ and the only Rees integer of  $bR$
is two. (In the notation of Discussion~\ref{3.1}, $D$ $=$ $E$ $=$
$V$, $J$ $=$ $N$, $m$ $=$ $2$, and  $R'[\theta]$ $=$ $R$.) However,
$R[P/b]'$ $=$ $R[\frac{y}{x}]'$ $\subseteq$ $k[[x,y]][\frac{y}{x}]$,
and  the powers of the maximal ideal  $M$ of $R$ define a valuation
on the quotient field of  $R$  that is readily seen to be a Rees
valuation ring of $P$,  but not a Rees valuation ring of $bR$ (since
$V$  is the only Rees valuation ring of $bR$). Therefore  $P$ $=$ $J
\cap R$ is not projectively equivalent to $bR$ $=$ $bV \cap R$, by
\cite[(3.4)]{CHRR}.

}
\end{exam}

With notation as in Example~\ref{failure}, the finite  integral
extension $A = R[x,y]$ contains an ideal $xA$ that is projectively
equivalent to $bA$ and the unique Rees integer of $xA$ is one. Thus
in relation to Question~\ref{QUES}, it seems natural to ask:

\begin{ques}
\label{prinx} {\em Let  $b$  be a nonzero nonunit in a Noetherian
integral domain  $R$. Does  there exist  a finite integral extension
domain  $A$ of $R$ having an ideal   $J$  whose Rees integers are
all equal to one such that $J$ is projectively equivalent to $bA$?

}
\end{ques}

With notation as in Discussion~\ref{3.1}, we give in
Proposition~\ref{equivs} several necessary and sufficient conditions
for the radical ideal $J \cap R'[\theta]$ to be projectively
equivalent to $bR'[\theta]$. The following definition is used in
this result.

\begin{defi}
\label{defiasymppd}
{\em
If  $I$  is a regular proper ideal in  $R$, then
${\hat A}^*(I)$  denotes  the set of
{\bf{asymptotic prime divisors}} of  $I$; that is,  ${\hat A}^*(I)$ $=$
$\{P \in \Spec(R) \mid P \in \Ass(R/(I^k)_a)$  for some positive
integer $k\}$.
}
\end{defi}

Concerning the hypothesis  ``$bR$ $=$ ${p_1}^{(m)} \cap \cdots \cap
{p_n}^{(m)}$'' in Proposition \ref{equivs}, it follows from either
Theorem \ref{prin}.1  or
Proposition \ref{multipleideals} below that, for each nonzero
nonunit $b$ in each Noetherian integral domain   $R$ there exists a
positive multiple $m$ of the Rees integers of  $bR$  and a
finite integral extension domain  $A_m$  of  $R$  such that
$b{A_m}'$ $=$ ${p_1}^{(m)} \cap \cdots \cap {p_n}^{(m)}$, where
$A_m'$ denotes the integral closure of $A_m$ and  $p_1,\dots,p_n$
are the prime divisors of $b{A_m}'$.

\begin{prop}
\label{equivs}
Let  $R$  be an integrally closed Noetherian domain,
let  $m$  be a positive integer,
let  $b$  be a nonzero nonunit in  $R$,
let  $p_1,\dots,p_n$  be
the (height one) prime divisors of  $bR$,
let  $J$ $=$ $p_1 \cap \cdots \cap p_n$, so  $J$
$=$ $\Rad(bR)$, and assume that  $bR$
$=$ ${p_1}^{(m)} \cap \cdots \cap {p_n}^{(m)}$.
Then the following statements are equivalent:

\noindent
{\bf{(\ref{equivs}.1)}}
$J$  is projectively equivalent to  $bR$.

\noindent
{\bf{(\ref{equivs}.2)}}
$J$  is projectively equivalent to some principal
ideal in  $R$.

\noindent
{\bf{(\ref{equivs}.3)}}
$bR$ $=$ $(({p_1} \cap \cdots \cap {p_n})^k)_a$  for some
positive integer  $k$.

\noindent
{\bf{(\ref{equivs}.4)}}
$(J^k)_a$  is principal for some positive integer  $k$.

\noindent
{\bf{(\ref{equivs}.5)}}
$J$  is invertible.

\noindent
{\bf{(\ref{equivs}.6)}}
$((p_1 \cap \cdots \cap p_n)^k)_a$ $=$
${p_1}^{(k)} \cap \cdots \cap {p_n}^{(k)}$  for all  positive integers  $k$.

\noindent
{\bf{(\ref{equivs}.7)}}
${\hat A}^*(J)$ $=$ $\{p_1,\dots,p_n\}$  (see (\ref{defiasymppd})).

\noindent
{\bf{(\ref{equivs}.8)}}
$J^kR_U \cap R$ $=$ $(J^k)_a$  for all positive integers  $k$,
where  $U$ $=$ $R - (p_1 \cup \cdots \cup p_n)$.
\end{prop}

\begin{proof}
Since  $J$ $=$ $p_1 \cap \cdots \cap p_n$  and
$bR$ $=$
${p_1}^{(m)} \cap \cdots \cap {p_n}^{(m)}$, by
hypothesis, it follows that if
(\ref{equivs}.6) holds, then  $(J^m)_a$ $=$
$bR$ $=$ $(bR)_a$  (since  $R$  is integrally
closed), so the case  $k$ $=$ $m$
of (\ref{equivs}.6)
implies that  (\ref{equivs}.1) holds.

It is clear
that (\ref{equivs}.1) $\Rightarrow$ (\ref{equivs}.2)

If (\ref{equivs}.2) holds, then let
$c$ $\in$ $R$  and let  $h,k$  be positive
integers  such that
$(J^k)_a$ $=$ $(c^hR)_a$, so  $(J^{k})_a$ $=$ $c^{h}R$,
since $R$  is integrally closed, so it follows
that  (\ref{equivs}.2) $\Rightarrow$ (\ref{equivs}.4).

Assume that (\ref{equivs}.4) holds and let
$c$ $\in$ $R$  such that
$(J^k)_a$ $=$ $cR$, so  $(J^{gk})_a$ $=$
$(c^gR)_a$ $=$ $c^{g}R$
for all positive integers  $g$, hence
$\Ass(R/(J^{gk})_a)$
$=$ $\Ass(R/c^{g}R)$
for all positive integers  $g$.
Since $R$  is integrally closed, it follows that
$\Ass(R/c^{g}R)$
is the set of height
one prime ideals in  $R$  that contain  $cR$,
so it follows that
$cR$ $\in$ $p_i$  for  $i$ $=$ $1,\dots,n$,
and  $p_1,\dots,p_n$  are the only height
one prime ideals in  $R$  that contain
$cR$, since
$\Ass(R/cR)$ $=$
$\Ass(R/(J^{k})_a)$.
Therefore, since $\Ass(R/{c^g}R)$ $=$
$\Ass(R/(J^{gk})_a)$
for all positive integers  $g$,
it follows from  Definition \ref{defiasymppd}
that (\ref{equivs}.7) holds, hence
(\ref{equivs}.4) $\Rightarrow$ (\ref{equivs}.7).

Assume that (\ref{equivs}.7) holds, so
$\Ass(R/(J^k)_a)$ $=$ $\{p_1,\dots,p_n\}$
for all positive integers  $k$, since each  $p_i$
is a minimal prime divisor of  $J$  and of
$(J^k)_a$.  Therefore for all positive
integers  $k$, $(J^k)_a$ $=$ $\cap\{(J^k)_aR_{p_i} \cap R
\mid i = 1,\dots,n\}$
$=$ $\cap\{J^kR_{p_i} \cap R
\mid i = 1,\dots,n\}$  (since  $I_a$ $=$ $I$
for all ideals in  $R_{p_i}$)
$=$ $\cap\{(J^kR_{p_i} \cap R_U) \cap R
\mid i = 1,\dots,n\}$
(where  $U$ $=$ $R - (p_1 \cup \cdots \cup p_n)$
$=$ $J^kR_U \cap R$, hence
(\ref{equivs}.7) $\Rightarrow$ (\ref{equivs}.8).

Let  $U$  be as in (\ref{equivs}.8).  Then it is
readily checked that
$J^kR_U \cap R$ $=$
${p_1}^{(k)} \cap \cdots \cap {p_n}^{(k)}$  for all  positive integers  $k$,
so  (\ref{equivs}.8) $\Rightarrow$ (\ref{equivs}.6)
(since  $J$ $=$ $p_1 \cap \cdots \cap p_n$).

The case $k$ $=$ $m$  of (\ref{equivs}.6) implies
that (\ref{equivs}.3) holds (with
$k$ $=$ $m$), since
$bR$ $=$
${p_1}^{(m)} \cap \cdots \cap {p_n}^{(m)}$, by
hypothesis, and
(\ref{equivs}.3) $\Rightarrow$ (\ref{equivs}.4),
since  $J$ $=$ $p_1 \cap \cdots \cap p_n$.

Finally, if (\ref{equivs}.2) holds, then  $J$  is
projectively equivalent to an invertible ideal,
so  $J$  is
invertible, by \cite[(2.10)(1)]{CHRR3}, so
(\ref{equivs}.2) $\Rightarrow$ (\ref{equivs}.5).
And if (\ref{equivs}.5) holds, then all ideals
that are projectively equivalent to  $J$  are
invertible, by \cite[(2.10)(1)]{CHRR3}, so  $\Ass(R/(J^k)_a)$
$=$ $\{p_1,\dots,p_n\}$  (the set of minimal prime
divisors of  $(J^k)_a$), by \cite[(3.9)]{CHRR3},
hence it follows from Definition
\ref{defiasymppd} that
(\ref{equivs}.5) $\Rightarrow$ (\ref{equivs}.7).
\end{proof}

\section{An application to asymptotic sequences.}
The main result in this section,
Proposition \ref{multipleideals},
partially extends Theorem \ref{prin}.1 to certain finite sets of
ideals, and its corollary
(\ref{coroasymp}) applies this to asymptotic sequences.
In the proofs we use the following definition.

\begin{defi}
\label{defiV} {\em Let  $I$  be a regular proper ideal in a
Noetherian ring  $R$ and let $k$  be a positive integer. Then the
{\bf{multiplicity}}  of  $k$  as a Rees integer of  $I$  is the
number of DVRs $(V,N)$ $\in$ $\Rees I$   such that $IV$ $=$ ${N}^k$.
}
\end{defi}

\begin{prop}
\label{multipleideals}
Let  $I_1,\dots,I_h$  be nonzero proper
ideals in a Noetherian domain  $R$ and, for  $i$ $=$ $1,\dots,h$,
let  $e_{i,1}, \dots, e_{i,n_i}$  be the Rees integers of  $I_i$
and  $m_i$ $=$ $e_{i,1} \cdots e_{i,n_i}$.
Assume that:  {\rm (a)}  $\Rees I_i$ $\cap$   $\Rees I_j$ $=$
$\varnothing$  for  $i$ $\ne$ $j$ in $\{1,\dots,h\}$. Then there
exists a simple free separable integral extension domain  $A$  of
$R$  such that, for  $i$ $=$ $1,\dots,h$, the Rees integers of $I_iA$
are all equal to  $m_i$.
\end{prop}

\begin{proof}
Let
$\Rees I_i$ $=$ $\{(V_{i,1},N_{i,1}),\dots,(V_{i,n_i},N_{i,n_i})\}$
(for  $i$ $=$ $1,\dots,h$), and let
$D$ $=$ $V_{1,1} \cap \cdots \cap V_{h,n_h}$, so
$D$  is a semi-local Principal Ideal Domain.  Also,
it follows from assumption (a)
that  $D$  has exactly
$n^*$ $=$ $n_1 + \cdots + n_h$  maximal ideals $M_{i,j}$
$=$ $N_{i,j} \cap D$, and  $D_{M_{i,j}}$ $=$ $V_{i,j}$  for
$i$ $=$ $1,\dots,h$  and  $j$ $=$ $1,\dots,n_i$.

By hypothesis, for  $i$ $=$ $1,\dots,h$ and  $j$ $=$
$1,\dots,n_i$, $I_iV_{i,j}$ $=$
${N_{i,j}}^{e_{i,j}}$, so
${M_{i,1}}^{e_{i,1}} \cap \cdots \cap {M_{i,n_i}}^{e_{i,n_i}}$
is an irredundant primary decomposition of  $I_iD$ (and
$e_{i,1},\dots,e_{i,n_i}$  are the Rees integers
of  $I_i$).  Let
$$ (\ref{multipleideals}.1) \qquad \qquad {e_{i,j}}^* =
\frac{m_i}{e_{i,j}}~~for~~i = 1,\dots,h~~and~~j=1,\dots,n_i$$
and let
$$ (\ref{multipleideals}.2) \quad \quad m = {e_{1,1}}^* \cdots
{e_{h,n_h}}^*,$$
so  $m$ $=$ ${m_1}^{n_1-1} \cdots {m_h}^{n_h-1}$.

Now resubscript the  $M_{i,j}$  as follows:
for  $i$ $=$ $1,\dots,h$  and  $j$ $=$
$1,\dots,n_i$  let
$M_{n_1 + \cdots + n_{i-1} + j}$ $=$ $M_{i,j}$
($M_{n_1 + \cdots + n_{i-1} + j}$ $=$ $M_{1,j}$,
if  $i$ $=$ $1$),
$e_{n_1 + \cdots + n_{h-1} + j}$
$=$ $e_{i,j}$, and ${e_{n_1 + \cdots + n_{h-1} + j}}^*$
$=$ ${e_{i,j}}^*$.  Then the remainder of the
proof of this proposition is similar to
the proof of Theorem \ref{prin}.2.
We  construct a chain
of semi-local Principal Ideal Domains
$$ (\ref{multipleideals}.3) \quad D = E_0 \subseteq E_1 \subseteq \cdots \subseteq
E_{n^*},$$
where, for  $k$ $=$ $1,\dots,n^*$, $E_k$  is the integral
closure of  $E_{k-1}$  in a realization  $L_k$  of
a realizable  ${e_k}^*$-consistent system  $S^{(k)}$ for $\{(E_{k-1})_N \mid
N$
is a maximal ideal in  $E_{k-1}\}$.
The systems  $S^{(k)}$  are all similar.  Specifically:
$S^{(k)}$  ramifies to the index  ${e_k}^*$  each of
the  ${e_1}^* \cdots {e_{k-1}}^*$ ($=$ $1$,
if  $k$ $=$ $1$) maximal ideals
$N$  in  $E_{k-1}$  that contract in  $D$  to $M_{k}$;
$S^{(k)}$  splits into  ${e_{k}}^*$  maximal ideals
all of the other maximal ideals
$N$  in  $E_{k-1}$; and, $S^{(k)}$  gives no
proper residue field extensions (that is, all of
the residue field extensions  $K_{i,j}$  (see
(1) in Definition \ref{consist})
of each maximal ideal  $N$  in  $E_{k-1}$
are chosen to be  $E_{k-1}/N$
$\cong$ $D/(N \cap D)$).

It is readily checked that each $S^{(k)}$  is
an ${e_k}^*$-consistent system for $\{(E_{k-1})_N \mid N$
is a maximal ideal in  $E_{k-1}\}$, and it is
realizable, by Theorem \ref{GK}(i).  Therefore
their ``composition'' yields the chain (\ref{multipleideals}.3)
of separable extensions of degrees
${e_1}^*,\dots,{e_{n^*}}^*$, resp., so the quotient
field  $L_{n^*}$  of  $E_{n^*}$
is separable over the quotient field
$L_0$  of  $R$  and  $D$, and  $[L_{n^*} : L_0]$
$=$ $m$  (with  $m$  as in (\ref{multipleideals}.2)).
It follows that each  ${M_{i,j}}$  is ramified to
the index  ${e_{i,j}}^*$  in each of the
$\frac{m}{{e_{i,j}}^*}$  maximal ideals in
$E_{n^*}$  that contain  $M_{i,j}$, so
${M_{i,j}}^{e_{i,j}} (E_{n^*})_N$ $=$
$(N^{{e_{i,j}}^*})^{e_{i,j}} (E_{n^*})_N$  for
each of these  $\frac{m}{{e_{i,j}}^*}$  maximal
ideals  $N$, and  ${e_{i,j}}^*e_{i,j}$ $=$ $m_i$, by
(\ref{multipleideals}.1).  Thus, for  $i$ $=$ $1,\dots,h$
and  $j$ $=$ $1,\dots,n_h$,
${M_{i,j}}^{e_{i,j}} E_{n^*}$  has only the
Rees integer  $m_i$  with multiplicity
$\frac{m}{{e_{i,j}}^*}$ (see Definition (\ref{defiV})).
Since  $I_i$ $=$ ${M_{i,1}}^{e_{i,1}} \cdots {M_{i,n_i}}^{e_{i,n_i}}$,
it follows that, for  $i$ $=$ $1,\dots,h$,
the Rees integers of  $I_iE_{n^*}$  are all equal to $m_i$.

Since $L_{n^*}$  is a finite
separable extension field of  $L_0$,
there exists an element  $\theta$  in
$L_{n^*}$
such that $L_{n^*}$ $=$ $L_0[\theta]$.
It is readily checked that this implies there exists
$r$ $\in$ $R$  such that  $r \theta$  is integral
over  $R$.  Therefore  $A$ $=$ $R[r \theta ]$  has
quotient field  $L_{n^*}$  and is a simple free
separable integral extension domain of  $R$.
Since the rings  $(E_{n^*})_N$  (with  $N$  a
maximal ideal in  $E_{n^*}$)  are the Rees
valuation rings of the ideals  $I_iA$,
by Remark \ref{fact}, it follows that,
for  $i$ $=$ $1,\dots,h$,
the Rees integers of  $I_iA$  are all equal to $m_i$.
\end{proof}

\begin{rema}
\label{rema1}
{\bf{(\ref{rema1}.1)}}
{\em
A similar proof shows that the
conclusion of Proposition \ref{multipleideals} continues to hold,
if the assumption (a)
is replaced with ``for all $i$ $\ne$ $j$ in $\{1,\dots,h\}$,
if  $(V,N)$ $\in$ $\Rees I_i$ $\cap$   $\Rees I_j$,
and if  $I_iV$ $=$ $N^{e_i}$
and $I_jV$ $=$ $N^{e_j}$, then  $m_i$  and  $m_j$  are chosen
so that  $e_jm_i$ $=$ $e_im_j$''.

\noindent
{\bf{(\ref{rema1}.2)}}
For  $i$ $=$ $1,\dots,h$  and  $j$ $=$ $1,\dots,n_i$ let
${e_{i,j}}^*$ $=$ $\frac{m_i}{e_{i,j}}$
(as in (\ref{multipleideals}.1)),
and let  $m$ $=$ ${e_{1,1}}^* \cdots {e_{h,n_h}}^*$
(as in (\ref{multipleideals}.2)).
Then it follows from the the proof of
Proposition \ref{multipleideals} that:
\newline
{\it{{\rm (*3)} the field  $L_{n^*}$ is a realization of the realizable}}
$m${\it{-consistent system}} $S^*$ $=$
$\{S_{1,1},\dots,S_{h,n_h}\}$, where, for $i$ $=$
$1,\dots,h$  and $j$ $=$ $1,\dots,n_i$,
$S_{i,j}$ $=$ $\{(K_{i,j,k},1,{e_{i,j}}^*) \mid k = 1,\dots,
\frac{m}{{e_{i,j}}^*}\}$.
}
\end{rema}

\begin{rema}
\label{alternaterema}
{\em
There is a simpler proof of Proposition \ref{multipleideals},
if at least one  $p_{i,j}$ $=$ $N_{i,j} \cap R$  is not
a maximal ideal.  Namely, in this case there exists an
algebraic extension field of  $D/M_{i,j}$  of
degree  $q$  for all integers  $q$.  (To see this,
there exists  $b$ $\in$ $I_{i,j}$
such that  $V_{i,j}$ $=$ ${C'}_p$, where  $C'$  is the
integral closure of  $C$ $=$ $R[I/b]$  and  $p$  is a
(height one) prime divisor of  $bC'$, so
$V_{i,j}/N_{i,j}$ $\cong$ $D/M_{i,j}$  is
a finite extension field of the quotient
field  $F$  of  $R/p_{i,j}$, by
\cite[(33.10)]{N2}.  Since  $R/p_{i,j}$  is
a Noetherian domain and not a field, there
exists a DVR with quotient field  $F$, so
there exists a DVR with quotient field  $D/M_{i,j}$,
and it is readily seen that  there exists an
algebraic extension field of  $D/M_{i,j}$  of
degree  $q$  for all integers  $q$.)  Therefore
let  $K_{i,j}$  be an algebraic extension field of
$D/M_{i,j}$  of degree  $\frac{m}{{e_{i,j}}^*}$
and let  $S_{i,j}$ $=$
$\{(K_{i,j},\frac{m}{{e_{i,j}}^*},{e_{i,j}}^*)\}$.
Also, for  $(i',j')$ $\ne$ $(i,j)$  let
$S_{i',j'}$ $=$
$\{(K_{i',j',k},1,{e_{i',j'}}^*) \mid k =
1,\dots,\frac{m}{{e_{i',j'}}^*}\}$.
Then it readily checked that:  (a)  $S$ $=$
$\{S_{1,1},\dots,S_{h,n_h}\}$  is an  $m$-consistent
system for  $\{(V_{1,1},\dots,V_{h,n_h}\}$; (b)  it
is realizable, by Theorem \ref{GK}(i);
(c)  $I_iE$ $=$ $(\Rad(I_iE))^{m_i}$  for  $i$ $=$
$1,\dots,h$, where  $E$  is the integral closure
of  $D$  in a realization  $L$  of  $S$; and,
(d) if  $A$ $=$ $R[r \theta]$  as in the last
paragraph of the proof of Proposition \ref{multipleideals},
then, for  $i$ $=$ $1,\dots,h$, the Rees valuation rings of
$I_iA$  are the rings in  $\{E_N \mid N$  is a maximal
ideal in  $E$  and  $I_i$ $\subseteq$ $N\}$, by
Remark \ref{fact}, so the Rees integers of  $I_iA$
are all equal to  $m_i$, by (c).
}
\end{rema}

To prove a corollary of Proposition \ref{multipleideals}, we
recall the following definition.

\begin{defi}
\label{defiasymp}
{\em Let  $R$  be a Noetherian ring  and let
$b_1,\dots,b_g$  be
regular nonunits in  $R$.
Then  $b_1,\dots,b_g$  are an
{\bf{asymptotic sequence in}}
$R$  in case
$(b_1,\dots,b_g)R$ $\ne$ $R$,
$b_1$  is not in any minimal prime
ideal in  $R$, and
$b_{i}$ $\notin$
$\cup \{P \mid P \in {\hat A}^*(b_1,\dots,b_{i-1})R\}$
for  $i$ $=$ $2,\dots,g$.
They are a
{\bf{permutable asymptotic sequence in}}
$R$  in case
each permutation of them is an
asymptotic sequence in  $R$.
}
\end{defi}

Concerning Definition \ref{defiasymp},
it is shown in \cite[(5.13)]{Mc} that every  $R$-sequence
is an asymptotic sequence, and
it is shown in \cite[(5.3)]{Mc} that if  $R$  is locally
quasi-unmixed, then an ideal is generated by
an asymptotic sequence if and only if it is an
ideal of the principal class.

\begin{coro}
\label{coroasymp}
Let  $b_1,\dots,b_g$  be an asymptotic sequence in a Noetherian
domain  $R$.  Then:

\noindent
{\bf{(\ref{coroasymp}.1)}}
For  $i$ $=$ $1,\dots,g$  let  $I_i$ $=$ $(b_1,\dots,b_i)R$
and let  $m_i$  be a positive common multiple of the
Rees integers of  $I_i$.
Then there exists a simple free separable
integral extension domain  $A$  of
$R$  such that, for  $i$ $=$ $1,\dots,g$, the Rees integers
of  $I_iA$  are all equal to  $m_i$.

\noindent
{\bf{(\ref{coroasymp}.2)}}
Assume that  $b_1,\dots,b_g$  is a permutable asymptotic sequence
and let  $\mathbf I$  be the set of all ideals of the form
$(b_{\pi(1)}, \dots, b_{\pi(k)})R$, where  $k$  varies over
$\{1,\dots,g\}$  and where  $\pi$  is an
arbitrary permutation of  $\{1,\dots,g\}$, so there are
$h$ $=$ $2^g-1$  ideals in  $\mathbf I$.  Let  $I_1,\dots,I_h$
be the ideals in  $\mathbf I$  and for  $i$ $=$ $1,\dots,h$
let  $m_i$  be a positive common multiple of the Rees
integers of  $I_i$.
Then there exists a simple free separable
integral extension domain  $A$  of
$R$  such that, for  $i$ $=$ $1,\dots,h$, the Rees integers
of  $I_iA$  are all equal to  $m_i$.
\end{coro}

\begin{proof}
For (\ref{coroasymp}.1),
$\Rees I_i$ $\cap$   $\Rees I_j$
$=$ $\varnothing$  for  $i$ $\ne$ $j$ in $\{1,\dots,g\}$,
since  $b_1,\dots,b_g$  is an asympototic sequence in  $R$,
so the conclusion follows from Proposition \ref{multipleideals}.

The proof of (\ref{coroasymp}.2) is similar, since
$b_1,\dots,b_g$  is a permutable asymptotic sequence in  $R$.
\end{proof}

\vspace{.15in}
\begin{flushleft}

Department of Mathematics, Purdue University, West Lafayette,
Indiana 47907-1395 {\em E-mail address: heinzer@math.purdue.edu}

\vspace{.15in}

Department of Mathematics, University of California, Riverside,
California 92521-0135
{\em E-mail address: ratliff@math.ucr.edu}

\vspace{.15in}

Department of Mathematics, University of California, Riverside,
California 92521-0135
{\em E-mail address: rush@math.ucr.edu}

\end{flushleft}

\end{document}